\documentstyle[12pt]{article}
\input{amssym.def}
\input{amssym.tex}

\newcommand{\es}{{\sf es}}
\def\newtheorems{\newtheorem{theorem}{Theorem}[section]

                 \newtheorem{lemma}[theorem]{Lemma}

                 }

\newtheorems

\begin{document}

\title{Martin's Axiom and $\Delta^2_1$ well-ordering of the reals}
\author{Uri Abraham\\ 
Department of Mathematics and Computer Science\\
Ben Gurion University, Be\'{e}r-Sheva, Israel;\\ 
Saharon Shelah\thanks{Partially sponsored by the Edmund Landau Center
for research in Mathematical Analysis, supported by the Minerva Foundation
(Germany), p.n. 458}\\
Institute of Mathematics\\
The Hebrew University, Jerusalem, Israel}
\maketitle
\begin{abstract} Assuming an inaccessible cardinal $\kappa$, there
is a generic extension in which $MA + 2^{\aleph_0} = \kappa$ holds 
and the reals
have a $\Delta^2_1$ well-ordering.\\
\end{abstract}

\section{Introduction}
The aim of this paper is to describe a technique that allows
the encoding of an arbitrary set of ordinals by a $\Delta^2_1$ formula in a
generic extension which is cofinality preserving.  This encoding is robust
enough to coexist with MA (Martin's Axiom).  Specifically, we will show, for
any model of ZFC set theory with an inaccessible cardinal $\kappa$, the 
existence of a cardinal
preserving generic extension in which $2^{\aleph_0} = \kappa + MA +$ {\em
there is a  $\;\Delta^2_1$ well-ordering of} ${\Bbb R}$.

Let us explain what is meant by a $\;\Delta^2_1$ well-order.  We refer here to
the structure $\langle H,\in \rangle$ where $H = H(\aleph_1)$ is the
collection of all hereditarily countable sets.  A $\Sigma^2_k$ formula is a
second-order formula of the form $\exists X_1 \subseteq H \ \forall X_2
\subseteq H \ldots \; \varphi(X_1,\ldots,X_k,\;a_1\ldots,a_n)$ with $k$
alternations of set quantifiers (unary predicates, $X_i$), and
where $\varphi$ is a first-order formula
(in which quantification is over $H$) with predicate names $X_1,\ldots,X_k$, and
variables $a_1,\ldots,a_n$ (which vary over $H$).  A $\;\Delta^2_k$ formula is
one that is equivalent to a $\Sigma^2_k$ and to a $\Pi^2_k$ formula.  A
 $\;\Delta^2_1$ well-ordering is one that is 
 given by a $\;\Delta^2_1$ formula $\psi(x,y)$
that defines a well-ordering
of ${\Bbb R}$.  Obviously, a $\Sigma^2_1$ linear ordering
of ${\Bbb R}$ is also a $\Pi^2_1$ ordering.

An alternative definition of $\Sigma^2_k$ formulas, which connects to the usual
definition of $\Sigma^1_n$ (projective) sets, is to look at third-order formulas
over $\langle {\Bbb N},+,\ldots\rangle$, that is, second order formulas over
${\Bbb R}$.

Our result cannot be improved to give a projective well-ordering of ${\Bbb R}$
because of a theorem of Shelah and Woodin \cite{ShWo}
which proves that there is no
well-ordering of ${\Bbb R}$ in $L({\Bbb R})$, assuming some large cardinal. 
Since any projective order is in $L[{\Bbb R}]$, and as a small extension, such
as the one described here, will not destroy any large cardinal above $\kappa$,
the $\Delta^2_1$ well-order cannot be improved to a projective well-order.

Though this paper can be read independently, 
it is a continuation of our \cite{AbSh}
work where another coding technique is described which does not add any new
reals.  Both that work and the present are motivated by a theorem of
Woodin \cite{Wo} which shows
 that if CH holds and there is a measurable cardinal
which is Woodin, then there is no $\Sigma^2_1$ well-order of the reals.  In
view of this result, a natural question is what happens if the CH is removed? 
Woodin has obtained the following result: Assuming an inaccessible cardinal
$\kappa$, there is a c.c.c. forcing extension in which $\kappa = 2^{\aleph_0}$
and
\begin{enumerate}
\item there is a $\Delta^2_1$ well-ordering of ${\Bbb R}$.
\item Martin's axiom holds for $\sigma$-centered posets.
\end{enumerate}

Since the poset used to get this extension has cardinality $\kappa$,
it does not destroy whatever large cardinal properties the ground model has
above $\kappa$, and hence the assumption of CH is necessary for Woodin's
theorem.

The theorem proved in this paper is a slight improvement of this theorem in
that MA replaces the restricted version for $\sigma$-centered posets,  but 
 our main point is to describe a different encoding technique.

We were also motivated by the following related result of Solovay:\\

There is a forcing poset of size $2^{2^{\aleph_0}}$ such that the following
holds in the extension
\begin{enumerate}
\item $2^{\aleph_0} = 2^{\aleph_1} = \aleph_2$.
\item $MA$ for $\sigma$-centered posets,
\item there is a $\Delta^2_1$ well-ordering of the reals.
\end{enumerate}

Let us emphasize that no inaccessible cardinal is needed for Solovay's result.
 Let us also mention here the main result of Abraham and Shelah \cite{AbSh}
\begin{quote}
There is a generic extension {\em that adds no new countable sets} in which
there exists a $\Sigma^2_2$ well-order of ${\Bbb R}$.\\
\end{quote}
The theorem proved in this paper will now be formally stated.\\
\noindent
{\bf Theorem}.  Let $\kappa$ be an inaccessible cardinal, and assume $GCH$
 holds below $\kappa$.  Then there is a forcing extension 
that changes no cofinalities and in which
\begin{enumerate}
\item $2^{\aleph_0} = \kappa +$ Martin's Axiom, and
\item there is a $\Delta^2_1$ well-ordering of ${\Bbb R}$.
\end{enumerate}

In a forthcoming work \cite{AbSh2} we will show that the inaccessible is
dispensable (but the continuum is $\aleph_2$ in this work).
\section{Overview}
  The idea of the proof is quite simple, and we first
 give  a general description. 
 The generic 
 extension is a length $\kappa$ mixed--support iteration consisting
of two components:  The first component iterates c.c.c. posets with the aim of
finally obtaining Martin's Axiom.  The second component is doing the coding.
  Quite arbitrarily, we have chosen
the set (called lim) 
of limit ordinals below $\kappa$ to be the support of the c.c.c.
component, and the set of successor ordinals (succ) to support the
coding component.  The iteration is a finite/Easton iteration.  This means
that the domain of each condition is finite on the limit ordinals, and has
cardinality $< \rho$ below every inaccessible cardinal $\rho \leq \kappa$.

For a regular cardinal $\lambda$, $F_\lambda$ denotes the club filter on
$\lambda$.  We say that a family $H\subseteq F_\lambda$ 
{\em generates} $F_\lambda$ iff  $\forall C \in F_\lambda \; 
\exists D \in H(D \subseteq C)$.  The least
cardinality of a generating family for $F_\lambda$ is called here ``the
generating number for $\lambda$''.  A crucial
question (in this paper) to ask about a
regular cardinal $\lambda$ is whether its generating number is $\lambda^{+}$ 
or higher:  it is through answers to these questions that the encoding works.

If $2^\lambda = \lambda^+$, then the generating number for 
 $\lambda$ is $\lambda^+$ of course, but it is easy to increase it by forcing, 
 say, $\lambda^{++}$ new subsets of $\lambda$ with conditions of size
$<\lambda$.  We denote with $C(\lambda,\mu)$ the poset that introduces $\mu$
subsets to $\lambda$ with conditions of size $<\lambda$.  $$C(\lambda,\mu) = 
\{f| \mbox{dom}(f) \subseteq \lambda \times \mu,\; \mbox{ range}(f) = 2, \mid 
f\mid < \lambda \}$$
where $\mid f \mid$ is the cardinality of the function $f$.
Equivalently, one can demand dom$(f) \subseteq \mu$ in the definition.
 Clearly $C(\lambda,\mu)$ is $\lambda$-closed, and if $\lambda^{<\lambda} =
\lambda$, then it satisfies the $\lambda^+ - $c.c.

The closure in $\lambda$ of each generic subset of $\lambda$ is a 
closed unbounded set that contains no old club set.
We will iterate such posets, varying $\lambda$, and taking care of MA as
well.

In the final generic extension, $2^{\aleph_0} = \kappa$, Martin's Axiom holds,
and the sequence of answers to the questions about the generating numbers for 
$\lambda < \kappa$ encodes a well-ordering
of ${\Bbb R}$ which is $\Delta^2_1$.
 As will be explained below, these questions are 
 asked only for even (infinite) successors
below $\kappa$, that is,
 cardinals of the form $\aleph_{\delta+2n}$ where $\delta
> 0$ is a limit ordinal and $1 \leq n < \omega$ (call this set of even
successor cardinals
\es).  It is convenient to use an
enumeration of \es\ that uses all the successor ordinals as indices: $\es =
\{ \lambda_j \mid j < \kappa$ {\it is a successor ordinal}$\}$.  So 
$\lambda_1 = \aleph_2$ is the first infinite even successor, $\lambda_2 =
\aleph_4,\ldots,\lambda_{\omega+1} = \aleph_{\omega +2}, \lambda_{\omega+2} =
\aleph_{\omega+4}$ etc.  In general, 
\begin{equation}
\label{1}
\mbox{if } \alpha = \delta + n+1 \ \mbox{where }\delta
\in \lim \ \mbox{and }n < \omega, \ \mbox{then }\lambda_\alpha
 = \aleph_{\delta + 2(n+1)}.
\end{equation}
 In the final model, the well-ordering of ${\Bbb R}$ is the sequence of
reals
 $\langle r_\xi | \xi < \kappa\rangle$ where $r_\xi \subseteq \omega$ is
 encoded  by setting $\alpha = \omega\xi$ and
\[ n \in r_\xi \;\mbox{iff the generating number for}\;\lambda =
 \lambda_{\alpha+n+1}\;\mbox{is}\;\lambda^{++} . \]

Why is it necessary to skip cardinals and to space the
$\lambda_\alpha$'s two cardinals apart?  Suppose that $r \subseteq \omega$ is
the first real we want to encode.  If $0 \in r$, then the first coding poset
is $C(\aleph_2,\aleph_4)$.  Recall that $GCH$ is assumed, and hence
cardinals are not collapsed, and $2^{\aleph_2} = \aleph_4$ after this forcing.
 Now if $1 \in r$, we may want to continue forcing 
 with $c(\aleph_3,\aleph_5)$, 
but this will collapse $\aleph_4$ since $2^{\aleph_2} = \aleph_4$.
  Thus we must start the next iteration at least two
cardinals apart, and forcing with $C(\aleph_4,\aleph_6)$ is fine.  
In general, $\lambda_{\alpha+1} = \lambda^{++}_\alpha$, enables
the proof that cardinals are not collapsed in the extension.

The coding component of our forcing will be an iteration of posets of type 
$C(\lambda_\alpha,\lambda_\alpha^{++})$ for well chosen $\alpha$'s.  This
choice will be made to obtain the desired coding by determining the
generating number for $\lambda \in \es$.

Let us take a closer, but still informal,
view of the forcing poset.
  If we denote  with $P_\alpha$ the
$\alpha$th stage of the iteration, then our final poset is $P_\kappa$.  For limit 
$\delta$'s, $P_\delta$ is the mixed support limit of $\langle P_i | i <
\delta\rangle$ with finite/Easton support.  This means that $f \in P_\delta$
iff $f$ is a partial function defined on $\delta$ such that $f\upharpoonright
i \in P_i$ for every $i < \delta$, and dom($f$) contains only finitely many
limit ordinals (this is the c.c.c. component), and $|\mbox{dom}(f) \cap \mu |
< \mu$ for any inaccessible cardinal $\mu$ (this is the Easton support
requirement of the coding component).  At successor stages $P_{j+1} \cong P_j
* Q_j$ is a two-step iteration, where $Q_j$ is a poset in $V^{P_j}$
characterized by the following.  For limit $j < \kappa,\; Q_j$ is in $V^{P_j}$
a c.c.c. forcing.  And for successor $j < \kappa$ of the form $\delta + i$, 
where $i \in \omega$ and $\delta \in \lim,\; Q_j$ is either the trivial poset,
or $C(\lambda_j,\lambda^{++}_j)$ which is the poset for adding
$\lambda^{++}_j$ many subsets to $\lambda_j = \aleph_{\delta +2i}$.  The
decision as to the character of $Q_j$ will be described later; the role of
$Q_j$ is to
encode one bit of information about some real.  
 This decision is made generically, 
 in $V^{P_j}$, and it depends on the real in $V^{P_j}$
that is being encoded.

So $P_1$ is some c.c.c. poset, and $P_2$ is $P_1$ followed by either the
trivial poset or by $C(\aleph_2,\aleph^{++}_2)$.  In the latter case, forcing
with $P_2$ makes $2^{\aleph_2} = \aleph_4$.  

The iteration continues in a similar fashion.  To illustrate one of the main
 points, let us see (only intuitively now) 
 why $\aleph_1$ is not collapsed.  We will show that every
$f : \omega_1 \to On$ in $V^{P_\kappa}$ (where $On$ is the class of
ordinals) has a countable approximation in
$V$, that is, a function $f'$ such that, for every $\alpha \in \omega_1$,
 $f(\alpha) \in f'(\alpha)$ where
$f'(\alpha)$ is a countable set of ordinals.

Observe first that the Easton component of $P_\kappa$ is $<\aleph_2$ closed. 
This means that if an increasing sequence $\langle p_i | i <
\omega_1\rangle$ of conditions in $P_\kappa$  have the same c.c.c.
component $(p_i \upharpoonright \lim = p_j \upharpoonright \lim$,
then there is an upper bound in
$P_\kappa$ to the sequence.  We say that $p$ is a {\em pure} extension of $q$
if $p$ extends $q$ and both have the same restriction to lim (same c.c.c.
component).  Now, if $f : \omega_1 \to On$ is a function in $V^{P_\kappa}$,
we define an increasing sequence $\langle p_i | i < \omega_1\rangle$ of
conditions in $P_\kappa$ such that $i < j \Rightarrow p_j$ is a pure extension
of $p_i$:  To obtain $p_{i+1}$ extend $p_i$ in countably many steps; at each
step find first an extension $q'$ of the previous step $q$ that forces a new
value for $f(i)$ (if there is one) and then take only the pure extension of
$q$ imposed by $q'$.  It turns out that this process will never take more than
countably many steps, or else we get a contradiction to the assumption that at
limit stages c.c.c. posets are iterated.  The upper bound $p \in P_\kappa$
of this pure increasing sequence ``knows'', for each $i < \omega_1$, all the
countable many possible values for $f(i)$.

We arrange the iteration in such a way that for every real $r \in
V^{P_\kappa}\;$ there is a unique 
limit ordinal $\delta = \delta(r)$ so that, for every $k \in \omega$, $k \in
r$ iff the generating number for $\lambda=\lambda_{\delta+(k+1)}$ is
$\lambda^{++}$ (by \ref{1}) $\lambda =
\aleph_{\delta + 2(k+1)})$. 

 Now the well-ordering
on ${\Bbb R}$ is defined by $$r_1 \prec r_2\ \mbox{iff}\ 
\delta(r_1) < \delta(r_2).$$  This
formula is certainly first-order expressible in $H(\kappa)$ (the collection
of sets o cardinality hereditarily $< \kappa$ in the
extension), but why is it $\Sigma^2_1$?  Why can we reduce it to
second--order quantification over $H9\aleph_1)$? The point is that 
$2^{\aleph_0}=\kappa \; + \; MA$, and we can speak
correctly within $H(\aleph_1)$ about $H(\kappa)$, and it takes a single
 second-order quantification to do that
(this trick was
used by Solovay in his theorem cited above; we will outline it now, and it
 will be explained in more
detail later.) 
 To express $r_1 \prec r_2$, just say:
\begin{quote} There is a relation $R$ over $H =
H(\aleph_1)$, such that $(H,R)$ satisfies enough of set theory (when $R$
interprets the membership relation $\in$), such that $R$ is well-founded and
such that every real is ``found''
 in $(H,R)$; moreover,  $(H,R)$ satisfies the following statement: 
``{\it every limit ordinal has the form $\delta(r)$ for some real $r$,
and $\delta( \hat r_1) < \delta(\hat r_2)$'',} where $\hat r$ is the
construction of $r \subseteq \omega$ in the model $(H,R)$.
\end{quote}
Since $R$ is well-founded, $(H,R)$ is collapsed to some $\in$ structure, $M$,
which turns out to be $H(\kappa)$ as we want.  The main points to notice in
order to prove this are
that (1) $M$ cannot contain less than $\kappa$ ordinals because it contains
all the reals, and a definable well-ordering of ${\Bbb R}$.
 (2) What $M$ considers to be a cardinal is really a cardinal,
because any possible collapsing 
 function in $H(\kappa)$ can be encoded by a real (with the almost
disjoint set technique which is applicable because of Martin's Axiom).  Since
this encoding real is in $M, H(\kappa)$ is included in $M$.  
 (3) $M$ does not
contain more ordinals than $\kappa$.  This is so since every limit ordinal
$\delta$ is connected to a single real which is encoded along the segment
$[\aleph_{\delta+2},\aleph_{\delta+\omega})$ by the characteristic of the club
filters.  Thus $M$ is $H(\kappa)$.

The details of this proof are written in the sequel.\\

\section{Mixed support iteration}

 In this section we describe how to iterate, with mixed support (Mitchell's
 type support),
c.c.c. posets and $\lambda$-complete posets, where  the
support of the c.c.c. component is finite, and the support of the complete
component is of Easton type---bounded below inaccessibles.

Let $\kappa$ be an inaccessible cardinal, and $\lambda < \kappa$ a
regular cardinal $> \aleph_1$. The non c.c.c posets in the iteration
are all assumed to be $\lambda$ closed. 
For definiteness we have chosen the support of
the c.c.c. posets to be the limit ordinals below $\kappa$, denoted here
lim (0 is in lim), and the $\lambda$-complete forcings are supported by the
successors below $\kappa$, denoted ``succ''.

For an ordinal $\mu \leq \kappa$, a {\em mixed support}\  iteration of
length $\mu$ is defined here to be a sequence of posets $\langle P_i | i \leq
\mu\rangle$ such that
\begin{enumerate}
\item The members of each $P_i$ are partial functions defined on $i$.
\item For limit $\delta \leq \mu,\;P_\delta$ is the mixed support {\em limit}
of $\langle P_i | i < \delta \rangle$.  This means the following.   $P_\delta$
consists of all the partial functions $f$ defined on $\delta$ such that
\begin{enumerate}
\item $f\upharpoonright i \in P_i$ for every $i < \delta$.
\item Dom$(f) \cap \lim$ is finite.
\item In case $\delta$ is inaccessible, $|\mbox{Dom}(f) \cap \mbox{succ}|<
\delta$.
\end{enumerate}
The partial order on $P_\delta$ is defined by $f \leq g$ iff for all
 $i < \delta\;f\upharpoonright i \leq g \upharpoonright i$ in $P_i$.
\item For successors $\eta + 1 \leq \mu,\;P_{\eta+1} \simeq P_\eta *
Q_{\eta}$ where $Q_\eta$ is a name of a poset in the universe of terms
$V^{P_\eta}$. So $f \in P_{\eta +1}$ iff $f\upharpoonright \eta \in P = P_\eta$
and $f\upharpoonright \eta \Vdash_P f(\eta) \in Q_{\eta}$.  The partial order
on $P_{\eta+1}$ is defined as usual.
\item For any limit ordinal $\delta < \mu,\;Q_\delta$ is in $V^{P_\delta}$ a
c.c.c. forcing (i.e., the empty condition in $P_\delta$ forces that).  For
successors $\alpha < \mu,\;Q_\alpha$ is $\lambda$-closed in $V^{P_\alpha}$
(closed under sequences of length $< \lambda$).
\end{enumerate}

The notation $\Vdash_\eta$ can be used for $\Vdash_{P_\eta}$.
It is convenient to define two conditions $p$ and $q$ in $P$ to be 
{\em
equivalent}\  iff they are compatible with the same conditions in $P$.  However,
it is customary not to deal with equivalence classes, and
 to write $p = q$ instead of $[p] = [q]$, and we shall accept
this convention.

For $i < \mu$ ($\mu$ is the length of the iteration) the restriction map $f
\mapsto f\upharpoonright i$ is a projection of $P_\mu$ onto $P_i$.  But for an
arbitrary set $A\subseteq i,\;f \upharpoonright A$ is not necessarily a
condition, and, even when it is a condition, it is possible that $[f] = [g]$ and
 $f \upharpoonright A \neq g \upharpoonright A$.  Therefore, the notation 
$f \upharpoonright A$
refers to the function $f$ itself and not to its equivalence class.

The set of functions $f \upharpoonright$ lim, for
$f \in P_\mu$, 
is called the ``c.c.c. component'' of $P_\mu$.  And the functions of the form 
$f \upharpoonright $ succ form the ``complete component'' of $P_\mu$. 
Let us
say that $f_2$ is a {\em pure} extension of $f_1$ in $P_\mu$ iff $f_1 \leq
f_2$ and $f_1  \upharpoonright \lim = f_2  \upharpoonright \lim$.  Thus, a
pure extension of $f_1$ does not touch the c.c.c. component.
(This definition refers to the functions $f_1$ and $f_2$ and not to their
equivalence classes in $P_\mu$.)

The following lemma is an obvious consequence of the assumed
$\lambda$-completeness of the posets in the complete component.\\

\begin{lemma}  $P_\mu$ is $< \lambda$ pure closed.  That is,
any purely increasing sequence $\langle
q_i | i < \tau\rangle$ of length $\tau < \lambda$ 
($q_j$ is a pure extension of $q_i$ for $i < j$) has a least upper bound in
$P_\mu$, which is a pure extension of each $q_i$.
\end{lemma}
Suppose now that $q \in P_\mu$, and $r$ is in the c.c.c. component of $P_\mu$.
 Then the sum $h = q + r$ is the function defined by

\[ h(i) = \left\{ \begin{array}{lcl}
r(i)&\;\mbox{if}\;& i \in \mbox{dom}(r) \\
q(i)&\;\mbox{if}\;& i \in \mbox{dom}(q)\setminus\mbox{dom}(r) \end{array}
\right. .\]
Whenever the notation $h=q + r$ is used, it is tacitly assumed that for every
$i,\;h \upharpoonright i \Vdash_i\; h(i) \in Q_i$ {\it and }
$r(i)$ {\it extends} $q(i)$. 
Hence $q + r \in P_\mu$ extends $q$. 
We have the following two easy lemmas on pure extensions given with no
proof.

\begin{lemma}  If $p_1 \leq p_2$ in $P_\mu$, then there is a pure extension $q$
of $p_1$ such that, setting $r = p_2 \upharpoonright \lim$, we have

\[ p_2 = q+r. \]
\end{lemma}
Thus any extension is a combination of a pure extension with a finitely
supported c.c.c. component.\\

\begin{lemma}  If $p_0 + r$ is a condition and $p_1$ is a pure extension
 of $p_0$, then $p_1 + r$ is a condition that extends $p_0 + r$.
\end{lemma}
The c.c.c. component of $P_\mu$ is certainly not a c.c.c. iteration, but the
following quasi c.c.c. property still carries over from the usual argument 
that iteration with finite support of c.c.c. posets is again c.c.c.\\

\begin{lemma}  Assume that $\omega_1$ is preserved by $P_{\mu'}$ for every 
$\mu' < \mu$.  Let $\{r_\xi | \xi < \omega_1\}$ be an uncountable subset of
the c.c.c. component of $P_\mu$.  If $q \in P_\mu$ is such that $q + r_\xi \in 
P_\mu$ can be formed for every $\xi < \omega_1$, then 
\begin{enumerate}
\item For some $\xi_1 \neq \xi_2,\;q + r_{\xi_1}$ and $q + r_{\xi_2}$ are
compatible in $P_\mu$.
\item There is some $r$ in the c.c.c. component of $P_\mu$ such that
 $q+r \in P_\mu$ and
\[
 \begin{array}{cl} q+r \Vdash_\mu & \mbox{there are unboundedly many}\; 
\xi < \omega_1 \\ & \mbox{with}\; q + r_\xi\in G\;\mbox{(the generic filter)}. 
\end{array}
\]
\end{enumerate}
\end{lemma}
{\bf Proof}.  Obviously, {\it (2)} implies {\it (1)} (because 
the posets are separative, and $p \Vdash  ``q + r_\xi \in G \, \mbox{''} $
implies $p_\xi \leq p$).
So we will only prove (2), by
induction on $\mu$. 

 Recall first that for any c.c.c. poset $Q$ and
uncountable subset 
$A \subseteq Q$ there is a condition $a \in A$ such that $a
\Vdash_Q\;A \cap G$ {\it is uncountable}. (Obvious warning: This does not
mean there are uncountably many $a' \in A$ 
with $a' \leq a$.)

If $\mu$ is limit, there is no problem in
using the familiar $\Delta$-argument in the case  $cf(\mu) = \omega_1$, and
the obvious application of the inductive assumption when $cf(\mu) \neq
\omega_1$.  For example, in case $cf(\mu) = \omega_1$, form a
$\Delta$-system out of dom$(r_\xi), \xi < \omega_1$, and let $d\subseteq i_0 <
\mu$ be the fixed finite core of the system.  Then apply the inductive
assumption to $q\upharpoonright i_0$ and to $r'_\xi = r_\xi \upharpoonright
d$, for $\xi$ in the $\Delta$ system. 
This gives some $r_0$ in the c.c.c. component of $P_{i_0}$ which satisfies
{\it 2} above for $q \upharpoonright i_0$ and the conditions $r'_\xi$. It is
not too difficult to see that $q+r_0$ is as required (use 
the fact that the c.c.c component of every condition has a finite support).

In case $\mu = j+1$ and $j$ is a limit ordinal (for this is the interesting
cse), then $P_\mu \simeq P_j * Q(j)$,
where $Q(j)$ is a c.c.c. poset in $V^{P_j}$.  Set $q' = q\upharpoonright j$,
and $r'_\xi = r_\xi \upharpoonright j$.  Apply induction to find $r'$ such
that 
$$ \begin{array}{cl}
q' + r' \Vdash_j & \mbox{\it for unboundedly many }\; 
\xi < \omega_1,\\ & q' + r'_\xi \in G_j\; \mbox{\it (the generic filter 
over }\; P_j). 
\end{array}
$$
 Then define a name $\sigma$ in $V^{P_j}$ of a subset of $\omega_1$ such
that

\[ [q' \Vdash_j \xi \in \sigma]\ \mbox{iff}\  q'+r'_\xi \in G_j . \]
Since \\
(1) $q' + r'$ forces that $\sigma$ is unbounded in $\omega_1$,\\
(2) $\omega_1$ is not collapsed in $V^{P_j}$ by our assumption,\\
(3) $Q(j)$ is c.c.c., \\
there is, by the remark made at the beginning of the proof,
 a name $a \in V^{P_j}$
such that $q'+r' \Vdash_j `` a$ {\it is some} $r_\xi(j) $ {\it for}
$r'_\xi$ {\it in} $G_j$ {\it such that} $a \Vdash_{Q(j)}$ (for unboundedly many
$ \zeta \in \sigma,\;r_\zeta(j) \in H$) ''. ($H$ is the $Q(j)$ generic filter.

Now it is immediate to combine $r'$ and $a$ to a function $r$ which is as
required.

The main property of the mixed support iteration is the following.
 \begin{lemma} Assume $P_\mu$ is a mixed support iteration as described
above of c.c.c. and $\lambda$-complete posets. For every cardinal
$\lambda' < \lambda$, every $f : \lambda'\to On$ in 
$V^{P_\mu}$  has a countable
approximation in $V$  (that is, a function $g$ defined on $\lambda'$ such that
for every $\alpha < \lambda',\;g(\alpha)$ is countable and $f(\alpha) \in
g(\alpha)$.)
\end{lemma}
{\bf Proof}.  By induction on $\mu$.  Observe first that the lemma implies
that any set of infinite cardinality $\lambda' < \lambda$ in the extension is 
covered by a ground model set of the same cardinality. Hence
cardinals $\leq \lambda$ are not collapsed in $V^{P_\mu}$. The lemma also
implies that,
for
regular uncountable $\lambda' < \lambda$, any club subset of $\lambda'$ in
$V^{P_\mu}$ contains an old club set in $V$.

  It is obvious that any c.c.c.
extension or $\lambda$-complete extension has the property described in the
theorem, namely
that functions on $\lambda'$ have countable approximations.  Hence,
in case $\mu = \mu_0 + 1$, the theorem is obvious:
First get the approximation in $V^{P_{\mu_0}}$ (assume without loss of
generality that the first approximation has the form $g : \lambda' \times
\aleph_0 \rightarrow On$, and then use induction to get a
second approximation in $V$. 

So assume that $\mu$ is a limit ordinal, and 
$f \in V^{P_\mu}$ is a function defined  on
$\lambda' < \lambda$.  We are going to define a pure increasing
sequence $\langle q_\xi | \xi < \lambda' \rangle$ in $P_\mu$ such that for
every $\alpha < \lambda'$ there is a countable set $g(\alpha)$ and

\[ q_{\alpha+1} \Vdash f(\alpha) \in g(\alpha) . \]
If this construction can be carried on, then  use the $< \lambda$ pure
completeness of $P_\mu$ to find an upper bound $q$ to this sequence.  Then 
$q \Vdash g$ {\it is a countable approximation to} $f$.

The definition of $q_{\xi +1}$ is done by defining 
 a pure increasing sequence 
$\langle q(\alpha) | \alpha < \alpha_0\rangle$ where $q(0) = q_\xi$, and 
 for
each $\alpha$, a finite function $r_\alpha$ in the c.c.c. component of $P_\mu$
 so that, for $\alpha \neq \alpha',\;q(\alpha) + r_\alpha$ and
$q(\alpha') + r_{\alpha'}$ force different values for $f(\xi)$.  The
definition of this sequence is continued
 as long as possible, and the following
argument shows that it must stop for some $\alpha_0 < \omega_1$, and then
$q_{\xi+1}$ is the pure supremum of this countable sequence,
and $g(\xi)$ is the set of all values forced there to be $f(\xi)$.
Indeed,
otherwise, $q(\alpha)$ can be defined for every $\alpha < \omega_1$ and we let
$q$ be the upper bound of this pure increasing sequence (recall that $\aleph_1
< \lambda$).  Then $q + r_\alpha$ is in $P_\mu$ for every $\alpha < \omega_1$
and it forces different values for $f(\xi)$.  This contradicts the quasi
c.c.c. lemma 3.4.

\section{Definition of the forcing extension}

 The description of the poset $P_\kappa$, used for the coding proof, is
given in this section by defining a
mixed-support iteration $\langle P_\mu |
\mu \leq \kappa\rangle$ as outlined in Section 2.

At successor stages: $P_{\mu+1} \cong P_\mu * Q_\mu$ where $Q_\mu$ is a poset
in $V^{P_\mu}$ defined thus.  If $\mu=\delta \in \lim$, then $Q_\delta$ is in
$V^{P_\delta}$ a c.c.c. poset of cardinality, say, $\leq \aleph_\delta$.  ($P_1$
 is a c.c.c. poset, say the countable Cohen poset.)  The choice of $Q_\delta$
is determined by some bookkeeping function, aimed to ensure that Martin's Axiom
holds in $V^{P_\kappa}$.
(The cardinality limitation is to ensure the right cardinalities to show
that cardinals are not collapsed.)

For successor ordinals of the form $j = \delta + i$ where $\delta$ is limit
and $0 < i < \omega,\;Q_j$ is defined to be in $V^{P_j}$ either the trivial
poset (containing a single condition) or the poset
$C(\lambda_j,\lambda_j^{++})$, where $\lambda_j = \aleph_{\delta+2i}$.  To
determine which alternative to take, define a function $g$ 
that gives, for every limit $\delta < \kappa$, a name $g(\delta) \in
 V^{P_\delta}$ such that, for every $\alpha < \kappa$, every real
 in
 $V^{P_\alpha}$ is some $g(\delta)$ for $\delta \geq \alpha$.
Suppose that $g(\delta)$ is interpreted as $r \subseteq \omega$ in
$V[G_\delta]$ (the generic extension via 
$P_\delta$); then this determines $Q_j$,
for every $j$ in the interval $(\delta , \delta + \omega)$, which has the form
$j = \delta + i_0 + 1$, by
\[ Q_j\;\mbox{is non-trivial iff}\ i_0 \in r . \]

In order to prove that $P_\kappa$ possesses the required properties (such as
not collapsing cardinals), we decompose $P_\kappa$ at any 
stage $\alpha < \kappa$, and write $P_\kappa \cong P_\alpha * P^\alpha_
\kappa$, where $P_\alpha$ is
the iteration up to $\alpha$, and $P^\alpha_\kappa$ is the remainder of the
iteration.  It is not hard to realize that $P^\alpha_\kappa$  is just like
 $P_\kappa$ except that $\lambda_1 = \aleph_2$ is replaced with 
 $\lambda_{\alpha
+1} = \aleph_{\alpha +2}$.  For this reason, we must first describe 
 $P^\alpha_\kappa$ and analyze its properties.

For each ordinal $\alpha < \kappa$ a mixed support iteration 
$\langle P^\alpha_\mu | \alpha \leq \mu \leq \kappa\rangle$ will be
defined by induction on
$\mu$.  The poset used to obtain the theorem is $P^0_\kappa$, but
the $P^\alpha_\kappa$ are necessary as well since the decomposition 
$P^0_\kappa \simeq P^0_\alpha * (P^\alpha_\kappa)^{V^{P^0_\alpha}}$ is used
to show the desirable properties of the iteration. 
This may also explain why we choose
the index $\mu$
of $P^\alpha_\mu$ to start from $\alpha$ and not from 0.  The
conditions in $P^\alpha_\mu$ are functions defined on the ordinal interval 
$[\alpha,\mu)$.

To begin with, $P^\alpha_\alpha$ is the trivial poset $\{\emptyset\}$
containing only one condition (the empty function).  The definition of $P^\alpha_{j+1} \simeq
 P^\alpha_j * Q^\alpha(j)$ depends on whether $j \in \lim$ or $j \in$ succ. 
If $j \in \lim$ then $Q^\alpha(j)$ is in $V^{P^\alpha_j}$ a c.c.c. poset of
cardinality $\leq \aleph_j$ (for definiteness).
 The choice of $Q^\alpha(j)$ for $j \in $ 
lim is determined by some bookkeeping
function which we do not specify now, the aim of which is to obtain Martin's
Axiom in $V^{P^0_\kappa}$.

If $j$ is a successor ordinal of the form $j =
\delta + i$ where $\delta$ is limit and $0 < i < \omega$, we require that
$Q^\alpha(j)$ is in $V^{P^\alpha_j}$ either the trivial poset, or 
$C(\lambda_j, \lambda^{++}_j)$ where $\lambda_j = \aleph_{\delta+2i}$ (all in
the sense of $V^{P^\alpha_j}$).  The exact description of $Q^\alpha(j)$
(i.e., the decision as to whether it is the trivial poset or the one that
introduces $\lambda^{++}_j$ club subsets to $\lambda_j$) is not 
needed to
 prove that cardinals are not
collapsed.

\begin{lemma}  For every successor $\alpha$, $P^\alpha_\mu$ is $\lambda_{\alpha}$ pure
closed. 
\end{lemma}
{\bf Proof}. The complete component of $P^\alpha_\mu$ consists of posets of
the form $C(\lambda_j,\lambda_j^{++})$ which are $\lambda_j$ closed.  Since 
$\lambda_\alpha \leq \lambda_j$ for all these $j$'s, the lemma follows.\\

\begin{lemma}
For every $\mu$ such that 
$\alpha \leq \mu \leq \kappa$,  $P^\alpha_\mu$ changes no
cofinalities and hence preserves cardinals.  In fact, this is deduced from the
following properties of the mixed support iteration $P_\mu^\alpha$.
\begin{enumerate}
\item For limit $\mu \leq \kappa$, the cardinality of $P^\alpha_\mu$ is $\leq
\aleph^+_\mu$, and if $\mu > \alpha$ is inaccessible, then $|P^\alpha_\mu| =
 \aleph_\mu$.
\item If $\mu$ is successor, $\mu = j+1$, then $P^\alpha_\mu$ satisfies the
$\lambda^+_j$-c.c. and its cardinality is $\leq \lambda^{++}_j$ (where
$\lambda_j = \aleph_{\delta+2i}$ if $j = \delta + i$ for $\delta$ limit and $0
\leq i < \omega$). Thus the GCH continues to hold in $V^{P^\alpha_\mu}$ for
$\lambda_j^+$ and higher cardinals.
\item For each $i$ such that
 $\alpha < i < \mu\;P^\alpha_\mu \cong P^\alpha_i * (P^i_\mu)^{V'}$ 
where $V'$ is $V^{P^\alpha_i}$.
\end{enumerate}
\end{lemma}
{\bf Proof}. Let us see first how 1,2,3 are used to show by induction 
that $P_\mu^\alpha$ preserves
cofinalities.  So let $g : \eta \to \sigma$ be a cofinal
function in $V^{P^\alpha_\mu}$ where $\eta$ is a regular cardinal.
We have to show that $cf(\sigma) \leq \eta$ in $V$ as well.
Assume first $\mu=j+1$ is a successor ordinal, and then 
$P^\alpha_\mu \cong P^\alpha_j * Q^\alpha(j)$.
 The case $j \in $ lim is obvious since $Q^\alpha(j)$ is then a c.c.c. poset. 
So assume that $j$ is a successor ordinal now, and $\lambda_j'$ is thus
defined.
The case $\lambda_j \leq \eta$ follows from the $\lambda_j^+$-c.c of $P_\mu^\alpha$. In case $\lambda_j > \eta$ use the
$\lambda_j$ completeness of $Q^\alpha(j)$ and induction. 

Now assume that $\mu$ is limit.
The proof divides into two cases.
Suppose, for some successor $j$ with $\alpha \leq j < \mu,\;\eta < \lambda_j$.  Then 
$P^\alpha_\mu \cong P^\alpha_j *
 (P^j_\mu)^{V'}$ where $V'$ is $V^{P^\alpha_j}$.  Lemma 3.5 was formulated
for quite a general mixed support iteration, and it can be applied in 
 $V'$ to $P^j_\mu$ 
to yield that the function $g$ has a countable approximation in $V'$. 
We may apply the inductive hypothesis and find an approximation of $g$ in
$V$.

In case $\eta \geq \lambda_j$ for all such $j$'s, $\eta \geq \aleph_\mu$. 
Apply cardinality or chain
condition arguments:  It follows in this case that $P^\alpha_\mu$ satisfies
the $\eta^+$-c.c. and hence the cofinality of $\sigma$ in $V$ is $\leq \eta$.

So now we prove the three properties
 by induction on $\mu$.  
The proof of {\it 1} and {\it 2} are fairly standard, and uses, besides the
definition of the Easton support, the inductive assumptions and the
restrictions on the cardinalities of the posets.

To prove {\it 3}, we shall define a map $f \mapsto \langle f\upharpoonright
i,\;f/i\,\rangle$ of $P^\alpha_\mu$ into $P^\alpha_i * (P^i_\mu)^{V'}$ as
follows.  Clearly, $f\upharpoonright i \in P^\alpha_i$.  To define the name 
$f/i$ in $V^{P^\alpha_i}$, we assume a $V$ generic filter, $G$, over
 $P^\alpha_i$, place ourselves in $V[G]$, and define
the function $(f/i)[G]$ which
interprets $f/i$ (for every $\xi \in
\mbox{dom}(f),\;f/i[G](\xi)$  is a name in $(P^i_\xi)^{V[G]}$ naturally
defined).  Let us check that this map is onto a dense subset of the
two-step iteration.  So let 
$\langle h,\tau\rangle \in P^\alpha_i * (P^i_\mu)^{V'}$.  By extending $h$ we
may assume that $h$ `knows' the finite domain of the c.c.c. component of
$\tau$.  That is, for some finite set $E_0 \subseteq
\mu,\; h \Vdash_i$ dom $(\tau)
 \cap\;\lim = E_0$.  Let  $E_1 = \{ \eta \in \mbox{succ}\;|\;$ some
extension of $h$ in $P^\alpha_i$ forces $\eta \in $ dom$(\tau)\}$.  Because
the cardinality of $P^\alpha_i$ is $< \aleph_{i+\omega},\;\;E_1$ is bounded
below inaccessible cardinals, and can serve as Easton support of a condition.
Now $f \in P_\mu^\alpha$ can be defined on $E_0 \cup E_1$, so that $\langle
h, f/i \, \rangle$ extends $\langle h, \tau \rangle$.

\section{The proof of the theorem}

 All the technical machinery is assembled, and we only have to
apply it.  The iteration has the form $P^0_\kappa$ and the definition of the
function $h$ that decides the value of $Q(j)$ is made so that  Martin's
Axiom holds in $V^{P^0_\kappa}$, and for every real $r$ in $V^{P^0_\kappa}$ 
there is a unique limit
ordinal $\delta(r)$ such that $$i \in r\; \mbox{iff for}\
j = \delta(r) + i + 1,\; Q(j)\; \mbox{is}\; C(\lambda_j,\lambda^{++}_j).$$

\begin{lemma}  For every successor $j < \kappa,\;Q(j)$ is
 $C(\lambda_j,\lambda^{++}_j)$ iff the club filter on $\lambda_j$ in
$V^{P^0_\kappa}$ has generating number $\lambda^{++}_j$ .\\
\end{lemma}
To prove the lemma, 
observe that any function $f : \lambda_j \to On$ has a countable
approximation in $P^0_{j+1}$.  This is so by Lemmas 4.1 and 3.5,
because $P^0_\kappa = 
P^0_{j+1} * (P^{j+1}_\kappa)^{V'}$, and $P_\kappa^{j+1}$ is in $V'$ a mixed
uspport iteration of c.c.c. and $\lambda_{j+1}$-closed posets. 
So every club subset of $\lambda_j$ in $V^{P^0_\kappa}$
contains a club in $P^0_j * Q(j)$, and then the generating number of
$\lambda_j$ in $V^{P^0_\kappa}$ and $V^{P^0_{j+1}}$ are the same. 
 But in $V^{P^0_j},\;2^{\lambda_j} = \lambda^+_j$ (by Lemma 4.2(2)),
 and hence the generating number in $V^{P^0_\kappa}$ is determined
in $P^0_{j} * Q(j)$ as follows.
If $Q(j)$ is trivial, then the generating number remains $\lambda_j^+$, and
if $Q(j)$ is $C(\lambda_j,\lambda_j^{++})$, then the generating number is
$\lambda_j^{++}$ of course.

The definition of the well-ordering of ${\Bbb R}$ in $V^{P^0_\kappa}$ is now
clear: $r_1 \prec r_2$ iff $\delta(r_1) < \delta(r_2)$.  Why is $\prec$ a
$\Sigma^2_1$ relation?  The answer was outlined in Section 2, and now 
more details are given.

 The ``almost disjoint sets encoding technique'' was introduced
by Solovay in \cite{So}, and the reader can find there a detailed exposition;
we only give an outline.
  Assume $\mu$ is a cardinal, and $s = \langle s_\xi |\xi <
\mu\rangle$ a collection of pairwise almost disjoint subsets of $\omega$.  Let
$X \subseteq \mu$ be any subset.  Then the following c.c.c. poset $P$ 
introduces a real $a \subseteq \omega$ such that, together with $s$, $a$
 encodes $X$. In fact, $\xi \in X$ iff $s_\xi \cap a$ is finite.

A condition $(e,c) \in P$ is a pair such that $e$ is a finite partial function
from $\omega$ to 2, and $c \subseteq X$ is finite.  The order relation
expresses the intuition that $e$ gives finite information on $a$, and $c$ is a
promise that for $\xi \in c$ the generic subset will not add any more members
of $a \cap s_\xi$.  So $(e_1,c_1)$
extends $(e_2,c_2)$ iff $e_2 \subseteq e_1,\;c_2 \subseteq c_1$, and for $\xi
\in c_2,\;s_\xi \cap E_1 \subseteq E_2$ (where $E_i = \{k|e_i(k) = 1\})$.

The intuitive meaning of this order relation becomes clear by the following
definition.  Let $G \subseteq P$ be generic; then set
\[ a = \{k | e(k) = 1\;\mbox{for some}\;(e,c) \in G\} . \]

It can be seen that, $a \cap s_\xi$ is finite for $\xi
\in X$, and is infinite for $\xi \not\in X$.

This almost disjoint set encoding is used to prove
that the $\Sigma^2_1$ definition given in Section 2 is really
equivalent to the well ordering 
$\prec$.  The main point is this.  Suppose Martin's Axiom
$+\, 2^{\aleph_0} = \kappa$, and $M$ is a transitive model of some part of ZFC
containing all the reals and a well-order of them (which is a class in $M$). 
Then $M$ contains all the bounded subsets of $\kappa$ as well.  Why? Well, let 
$X \subseteq \mu < \kappa$ be any bounded set.  Since $M$ contains a set of
$\mu$ reals, it also contains a sequence of $\mu$ pairwise almost disjoint
subsets of $\omega$ (taken, for example, as branches of
$2^{\stackrel{\omega}{\smile}}$).  By Martin's Axiom, there is a set $a
\subseteq \omega$ that encodes $X$.  As $a \in M,\;X \in M$ as well.

\section{A weakening of the GCH assumption}
The theorem required GCH (below $\kappa$)
to ensure that cardinal are not collapsed. In this section this assumption is 
weakened somewhat in demanding that $2^\mu=\mu^+$ only on some closed
unbounded set of cardinals $\mu < \kappa$.

To see this, let $\langle \mu_i \mid i < \kappa \rangle$ be an enumeration
of a club set of limit cardinals, such that $2^{\mu_i}=\mu_i^+$, and
cf$(\mu_{i+1}) > \mu_i^+$, and $(\mu_{i+1}) ^{\leq \mu_i}= \mu_{i+1}$.

The construction is basically the same as before, but
 $\mu_i$ replaces $\lambda_i$ 
and  the main point is this:
For a successor $j=\delta+i$, where $\delta < \kappa$ is limit
and $0< i < \omega$, $P_{i+1}=P_i * Q(j)$ where $Q(j)$ is now a poset that
adds either $\mu_{j}$ or $\mu_{j+1}$ 
subsets to $\mu_j^+$. Now if $M$ is as before a transitive
model that contains all the reals, then the club sequence can be reconstructed
by asking the questions about the generating numbers. If one starts with
$\mu_0$, then the original sequence is reconstructed; starting with another
cardinal may result in another club. However, this club intersects the
original sequence of the $\mu_i$'s, and hence both sequences have an
equal end-section. Hence we must demand that the well-ordering of
$\Bbb R$ is determined by any end section of the club.

\end{document}